\documentclass[11pt,reqno,tbtags,a4paper]{amsart}
\usepackage{amssymb}
\usepackage{url}
\usepackage[square,numbers]{natbib}
\bibpunct[, ]{[}{]}{;}{n}{,}{,}

\title%[]
{On convergence for graphexes}

\date{21 February, 2017; revised 23 March, 2022}
\author{Svante Janson}
\address{Department of Mathematics, Uppsala University, PO Box 480,
SE-751~06 Uppsala, Sweden}
\email{svante.janson@math.uu.se}
\urladdr{http://www.math.uu.se/svante-janson}

\subjclass[2010]{} 
\numberwithin{equation}{section}

\renewcommand\le{\leqslant}
\renewcommand\ge{\geqslant}

\allowdisplaybreaks

\theoremstyle{plain}% default
\newtheorem{theorem}{Theorem}[section]
\newtheorem{lemma}[theorem]{Lemma}

\theoremstyle{definition}

\newtheorem{problem}[theorem]{Problem}
\newtheorem{remark}[theorem]{Remark}

\theoremstyle{remark}

\newenvironment{romenumerate}[1][-10pt]{% optional argument changes indentation
\addtolength{\leftmargini}{#1}\begin{enumerate}% gives (i), (ii) etc.
 \renewcommand{\labelenumi}{\textup{(\roman{enumi})}}%
 \renewcommand{\theenumi}{\textup{(\roman{enumi})}}%
 }{\end{enumerate}}

\newcounter{oldenumi}
{\setcounter{oldenumi}{\value{enumi}}
\begin{romenumerate} \setcounter{enumi}{\value{oldenumi}}}
{\end{romenumerate}}

\newcounter{thmenumerate}

\newcounter{xenumerate}   %no left indentation; thus wider lines

\newcommand{\refT}[1]{Theorem~\ref{#1}}

\newcommand{\refL}[1]{Lemma~\ref{#1}}

\newcommand{\refS}[1]{Section~\ref{#1}}
\newcommand{\refSS}[1]{Section~\ref{#1}}
\newcommand\REM[1]{{\raggedright\texttt{[#1]}\par\marginal{XXX}}}
\newcommand\XREM[1]{\relax}

\begingroup
  \divide\count255 by 60
  \multiply\count255 by -60
  \advance\count255 by \time
  \ifnum \count255 < 10 \xdef\klockan{\the\count1.0\the\count255}
  \else\xdef\klockan{\the\count1.\the\count255}\fi
\endgroup

\newcommand{\summo}{\sum_{m=0}^\infty}

\newcommand\set[1]{\ensuremath{\{#1\}}}
\newcommand\bigset[1]{\ensuremath{\bigl\{#1\bigr\}}}
\newcommand\Bigset[1]{\ensuremath{\Bigl\{#1\Bigr\}}}

\newcommand\bigpar[1]{\bigl(#1\bigr)}

\newcommand\xcpar[1]{\{#1\}}

\newcommand\bigabs[1]{\bigl|#1\bigr|}

\newcommand\lrabs[1]{\left|#1\right|}
\def\rompar(#1){\textup(#1\textup)}    % usage: \rompar(...)

\def\xexp(#1){e^{#1}}

\newcommand\ntoo{\ensuremath{{n\to\infty}}}

\newcommand\stoo{\ensuremath{{s\to\infty}}}
\newcommand\ttoo{\ensuremath{{t\to\infty}}}

\newcommand\norm[1]{\|#1\|}

\newcommand\punkt{.\spacefactor=1000}    % om problem!
\newcommand\iid{i.i.d\punkt}    
\newcommand\ie{i.e\punkt}
\newcommand\eg{e.g\punkt}

\newcommand\cf{cf\punkt}
\newcommand{\as}{a.s\punkt}
\newcommand{\aex}{a.e\punkt}
\newcommand{\tend}{\longrightarrow}
\newcommand\dto{\overset{\mathrm{d}}{\tend}}

\newcommand\asto{\overset{\mathrm{a.s.}}{\tend}}
\newcommand\eqd{\overset{\mathrm{d}}{=}}

\newcommand\bbR{\mathbb R}

\newcommand\bbN{\mathbb N}

\newcommand\bbQ{\mathbb Q}
\newcommand\bbZ{\mathbb Z}

\newcounter{CC}
\newcounter{cc}
\newcommand\E{\operatorname{\mathbb E{}}}
\renewcommand\P{\operatorname{\mathbb P{}}}

\newcommand\Po{\operatorname{Po}}

\newcommand\Be{\operatorname{Be}}

\newcommand\gd{\delta}
\newcommand\gD{\Delta}
\newcommand\gf{\varphi}

\newcommand\gG{\Gamma}

\newcommand\gl{\lambda}

\newcommand\gs{\sigma}

\newcommand\gth{\theta}
\newcommand\eps{\varepsilon}

\renewcommand\phi{\xxx}  %% WARNING

\newcommand\cA{\mathcal A}

\newcommand\cF{\mathcal F}
\newcommand\cG{\mathcal G}

\newcommand\cL{{\mathcal L}}
\newcommand\cM{\mathcal M}
\newcommand\cN{\mathcal N}

\newcommand\cP{\mathcal P}

\newcommand\cS{{\mathcal S}}

\newcommand\cW{\mathcal W}

\newcommand\ett[1]{\boldsymbol1\xcpar{#1}} 
\newcommand\indic[1]{\boldsymbol1\xcpar{#1}}

\newcommand\etta{\boldsymbol1}

\newcommand\qw{^{-1}}

\newcommand\qqw{^{-1/2}}

\newcommand\intoo{\int_0^\infty}

\newcommand\oi{\ensuremath{[0,1]}}

\newcommand\oii{[0,1]^2}
\newcommand\ooo{[0,\infty)}

\newcommand\setoi{\set{0,1}}

\newcommand\dd{\,\mathrm{d}}

\newcommand\lhs{left-hand side}
\newcommand\rhs{right-hand side}

\newcommand\bbQp{\bbQ_+}
\newcommand\bbZp{\bbZ_+}
\newcommand\bbRp{\bbR_+}
\newcommand\bbRpp{\bbRp^2}
\newcommand\restr[1]{|_{#1}}
\newcommand\restrr{\restr{r}}
\newcommand\lbl{\mathsf{Lbl}}
\newcommand\MRR{\cM(\bbRp^2)}
\newcommand\PMRR{\cP(\MRR)}

\newcommand\cNN{\cN_{\mathsf{s}}}
\newcommand\cNNs{\cN_{\mathsf{s},\mathsf{s}}}
\newcommand\Adj{\cNNs(\bbRpp)}
\newcommand\vth{\vartheta}
\newcommand\fG{\mathfrak{G}}
\newcommand\fGf{\mathfrak{G}_\mathsf{f}}
\newcommand\fGl{\widehat{\mathfrak{G}}}
\newcommand\WW{\mathfrak{W}}
\newcommand\WWX{\widetilde{\WW}}
\newcommand\dgp{d_{\mathsf{GP}}}
\newcommand\dgs{d_{\mathsf{GS}}}
\newcommand\togp{\to_{\mathsf{GP}}}
\newcommand\togs{\to_{\mathsf{GS}}}
\newcommand\chS{\hat{\cS}}
\newcommand\os{[0,s]}
\newcommand\oss{[0,s]^2}
\newcommand\lbls{\lbl_s}
\newcommand\lblr{\lbl_r}
\newcommand\rr{^{(r)}}
\newcommand\hW{\hat W}
\newcommand\rs{_r^{(s)}}
\newcommand\hWss{\hW_{G_s,s}}
\newcommand\bX{\bar X}
\newcommand\bM{\bar M}
\newcommand\GG{\overline{G}}
\newcommand\tocut{\to_{\dcut}}
\newcommand\tocuts{\to_{\dcuts}}
\newcommand\cn[1]{\norm{#1}\cut}
\newcommand\cut{_{\square}}
\newcommand\dcut{\delta_{\square}}
\newcommand\dcuts{\delta_{\square}^{\mathsf s}}
\newcommand\mpp{measure-preserving}

\newcommand\mpb{\mpp{} bijection}
\newcommand\NN{^{(N)}}
\newcommand\str[1]{^{(#1)}}
\newcommand\strc{\str c}
\newcommand{\Lovasz}{Lov\'asz}

\begin{document}

\begin{abstract} 
We study four different notions of convergence for graphexes,
recently introduced by Borgs, Chayes, Cohn and Holden, and by Veitch and Roy.
We give some properties of them and some relations between them. We also
extend results by Veitch and Roy on convergence of empirical graphons.
\end{abstract}

\maketitle

\section{Introduction}\label{S:intro}

The theory of graph limits for dense graphs and the representation of such graph
limits by classical graphons has developed over the last decade and has been
very succesful, see \eg{} the book by \citet{Lovasz}.
Here, a classical graphon is a symmetric measurable function
$W:S^2\to\oi$, where $S$ is a probability space; without loss of generality
one can take $S=\oi$.

There have been many different attempts to find corresponding results for
sparse graphs. One recent approach has been through random graphs defined by 
certain discrete exchangeable random measures on $\bbRpp$.
Exchangeable random measures on $\bbRpp$ were characterized by
\citet{Kallenberg1990,Kallenberg2}, and his general construction of such
measures can be interpreted as a construction of random graphs.
(Note that the classical theory of graph limits and graphons on $\oi$ 
can be derived from the related characterizations by Aldous and Hoover of
exchangeable arrays, see \cite{Austin,SJ209}, although this was not
the original method or motivation.)
This use of Kallenberg's construction 
was first done by \citet{CF} in a special case, and extended by 
\citet{HSM15}, \citet{VR} and \citet{BCCH16}.
We will here use the 
version by  \citet{VR},  described in detail in \refS{SSgraphon}
below. It uses a \emph{graphex}, which is a triple $(I,S,W)$, where the most
interesting part is $W$ which is a graphon, but a graphon in a new more
general sense; $W$ is defined on the infinite measure space
$\bbRp$ instead of $\oi$.

More generally, graphons can be defined on  
any $\gs$-finite measure space $S$. This was developed by
\citet{BCCH16}.
However, one of their results is that it is possible to take $S=\bbRp$
without loss of generality, and we will in the present paper only consider
this case (following \cite{VR}).

Having defined graphons and graphexes, it is natural to define a topology on
them, and thus a notion of convergence. 
For  classical  graphons, there are several quite different ways to define
convergence, but they are all equivalent,
see \eg{} \citet{BCLSV1,BCLSV2}.
(This important fact is closely related to the fact that the space of
classical graphons, modulo equivalence, is compact.)

In the present, more general, context, there are also several possibilities,
but, unforunately, they are not equivalent.
It seems not yet clear which notion(s) of convergence that will turn out to be
useful in applications, and it seems that several possibilities ought to be
studied more. The present paper is a small contribution to this.

We consider in this paper four different notions of convergence for
graphexes and graphons.
Two of them, denoted  $\togp$ and $\togs$, were defined by
\citet{VR2}, based on convergence in distribution of the corresponding random
graphs; we stress 
(which is implicit in \cite{VR2})
that both convergences are metric, \ie, can be defined
by (pseudo)metrics. 
The two other notions of convergence apply only to
the special case of integrable graphons;
they use the (pseudo)metrics $\dcut$ and $\dcuts$ defined by \citet{BCCH16}
(see also \cite{SJ311}).
See \refS{Sconv} for detailed definitions.
(The four metrics studied in the present are not the only possible ones.
In particular, we do not consider the left convergence studied in
\cite{BCCH16}.)

We show that for integrable graphons, convergence in $\dcut$ ($\dcuts$)
implies convergence $\togp$ ($\togs$). We conjecture that the converses do
not hold, but we leave that as an open problem.

Each graphex defines a random graph process $(G_s(\cW))_{s\ge0}$, see
\refS{SSgraphon} below.
\citet[Theorem 2.23]{BCCH16} show that for any integrable graphon $W$, the
empirical graphon defined by the random graph $G_s(W)$ \as{} converges to 
$W$ in the metric $\dcuts$ as \stoo. Similarly, 
\citet{VR2} 
show that for any graphex $\cW$, 
the empirical graphon defined by $G_s(\cW)$ converges to the graphex $\cW$ 
in $\togp$ (after suitable stretching) and in $\togs$
as
\stoo; however, they prove this only for a sequence $s_k\to\infty$, and in
general only with convergence of probability.
We extend their theorems to convergence for the full family $G_s$
with a continuous parameter $\stoo$; moreover, we show \as{} convergence.
(See Theorems \ref{T2} and \ref{T3}.)
Furthermore, in order to prove this result,
\citet{VR2} first show a related convergence result for a randomly relabelled
version of the random graph, again for sequences $s_k\to\infty$.
Again we improve their result to
convergence for the continuous parameter $\stoo$
(\refT{T1}).
(In both cases, we make only some minor technical improvements in the proof;
the proofs are thus essentially due to \cite{VR2}.)

\section{Notation and preliminaries}\label{Snot}

Much of the notation follows \citet{VR,VR2}, but there are various
modifications and additions for our purposes.

$\gl$ denotes Lebesgue measure (in one or several dimensions).

$\bbRp:=\ooo$, the set of non-negative real numbers.

If $S$ is a measurable space, then $\cP(S)$ is the set of probability
measures on $S$.
If $X$ is a random variable in some measurable space $S$, then
$\cL(X)$ denotes the distribution of $X$; thus $\cL(X)\in\cP(S)$.
If $S$ is a metric (or metrizable) space; we equip $\cP(S)$ with the usual
weak topology,
see \eg{} \cite{Billingsley} or
\cite[Chapter 4]{Kallenberg}.
Note that if $S$ is a Polish space (\ie, it can be given a complete and
separable metric), then $\cP(S)$ is Polish too,
see \cite[Appendix III]{Billingsley}.
We denote  convergence in distribution of random variables in $S$ by $\dto$;
recall that $X_n\dto X$ means $\cL(X_n)\to\cL(X)$ in $\cP(S)$.

If $X$ and $Y$ are random variables (defined on the same probability space
and, for simplicity, with values in some Polish spaces $\cS_X$ and $\cS_Y$),
then $\cL(X\mid Y)$ denotes the conditional distribution of $X$ given $Y$;
note that this is a random probability measure on $\cS_X$ that can be
regarded as a function of $Y$. We use also $(X\mid Y)$ for a random variable
with this conditional distribution.

If furthermore $S$ is a locally compact Polish space
(= locally compact second countable Hausdorff space),
then
$\cM(S)$ is the set of locally finite Borel measures on
$S$. (We will only use $S=\bbRpp$ and subsets thereof.) 
We equip $\cM(S)$ with the vague topology, 
which makes $\cM(S)$ into a Polish space,
see
\cite[Appendix A.2 and Theorem A.2.3]{Kallenberg}.
Furthermore,  $\cN(S)$ is the subset of integer-valued
measures in $\cM(S)$, 
\ie, the set of all locally finite sums of unit point masses
$\gd_{x}$, and we let $\cNN(S)$  be the subset of simple
integer-valued measures, \ie, locally finite sums of distinct unit point masses.
It is easily seen that $\cN(S)$ and $\cNN(S)$ are measurable subsets of
$\cM(S)$.

\subsection{Graphs and adjacency measures}

We consider both unlabelled and labelled graphs; in the labelled case, each
vertex is labelled with a real number in $\bbRp$, and these labels are
supposed to be distinct.
The graphs may be finite or countably infinite, but we always assume that
there are no isolated vertices; thus a graph $G$ is determined by its edge
set $E(G)$. We furthermore consider only graphs that are simple in the sense
that there are no multiple edges; in general we allow loops 
(but in many applications we do not have any).

We denote the vertex set and edge set of a graph $G$ by $V(G)$ and $E(G)$,
and let $v(G):=|V(G)|$ and $e(G):=|E(G)|$ be the numbers of vertices and edges.

If $\gG$ is a labelled graph, 
then the corresponding unlabelled graph, obtained by
ignoring the labels, is denoted $\cG(\gG)$. 
Conversely, if $G$ is an unlabelled graph
and $s>0$, 
then $\lbl_s(G)$ is the (random) labelled graph obtained by labelling the
vertices by 
random \iid{} labels that are $U(0,s)$, \ie, uniformly distributed in $(0,s)$.
(Note that this yields distinct labels a.s., so we may assume that the
labels are distinct as required above.)
If $G$ is a labelled graph, we define $\lbl_s(G)$ in the same way; thus
relabelling the vertices randomly (regardless of their original labels).
In other words,  $\lbl_s(G):=\lbl_s(\cG(G))$.
When $G$ is a random graph (labelled or not), $\lbl_s(G)$ is defined by
taking the labelling independent of $G$.

If $\gG$ is a labelled graph, 
we represent the edge set $E(\gG)$ of $\gG$ (and thus the graph $\gG$ itself)
by the measure 
\begin{equation}\label{xi}
\xi=\xi(\gG):=\sum_{x,y\in\bbRp:(x,y)\in E(\gG)}\gd_{(x,y)}  
\end{equation}
on $\bbRp^2$; an edge between two distinct vertices 
labelled $x$ and $y$ is thus represented by the two point masses
$\gd_{(x,y)}+\gd_{(y,x)}$, while a loop (if such exist) at a vertex
labelled $x$ is represented by $\gd_{(x,x)}$. (We consider undirected
graphs, and thus the endpoints of an edge have to be treated symmetrically.)
Note that $\gG$ is determined (as a labelled graph) by $\xi$.

If $\gG$ is a labelled graph, and $\xi$ the corresponding measure, then
for $r\ge0$, $\gG\restrr$ denotes the induced subgraph of $\gG$ obtained by
first eliminating all vertices with labels $>r$, and any edges incident to
such a vertex, and then also removing all remaining vertices that have become
isolated. In other words, we keep the edges whose endpoints both have labels
$\le r$, and the endpoints of these edges. We let $\xi\restrr$ denote the
corresponding measure on $\bbRp^2$, and note that this is just the
restriction of $\xi$ to $[0,r]^2$. 

A labelled graph $\gG$ is \emph{locally finite} if $\gG\restrr$ is finite
for each $r<\infty$;  equivalently, $\xi(\gG)$ is a locally finite measure.
We consider only locally finite graphs.
We say that a measure $\xi$ on $\bbRpp$ is an \emph{adjacency measure} if it
is given by \eqref{xi} for some locally finite labelled graph $\gG$.
Hence, a measure is an adjacency measure if and only if it is a symmetric
measure in $\cNN(\bbRpp)$. We denote the set of adjacency measures by
$\Adj$, and let $\fGl$ be the set of  locally finite labelled graphs.

Thus, \eqref{xi} defines a 1--1 correspondence $\gG\leftrightarrow\xi(\gG)$
between the sets $\fGl$ and $\Adj$ of locally finite labelled graphs
and adjacency measures.
We give the set of adjacency measures $\Adj$ the subspace topology as a
subset of $\cM(\bbRpp)$, and give $\fGl$ the corresponding topology induced
by the correspondence \eqref{xi}.
Thus $\Adj$ and $\fGl$ are metric spaces, and a sequence $\gG_n\to\gG$ in
$\fGl$ if and only if $\xi(\gG_n)\to\xi(\gG)$ in $\cM(\bbRpp)$, \ie, 
in the vague topology.
%To emphazise this, we will use the notation $\gG_n\dto\gG$ in $\MRR$

Returning to unlabelled graphs, we let $\fG$
be the set of finite or countably infinite unlabelled graphs,
and 
$\fGf$ the subset of finite unlabelled graphs.
Then $\fGf$ is countable, and we give it the discrete topology.
(We do not define a topology on $\fG$.)

\subsection{Graphons, graphexes and random graphs}\label{SSgraphon}

A \emph{graphon} is (in the present context) a 
symmetric, measurable function $W:\bbRp^2\to\oi$ that satisfies the
integrability conditions, where $\mu_W(x):=\intoo W(x,y)\dd y$:
\begin{romenumerate}
\item \label{graphona}
$\mu_W(x)<\infty$ for \aex{} $x$
and
$\gl\set{x:\mu_W(x)>1}<\infty$;
\item 
$\int_{\bbRp^2}W(x,y)\indic{\mu_W(x)\le 1}\indic{\mu_W(y)\le 1}\dd x\dd
y<\infty;
$
\item \label{graphono}
$\int_{\bbRp} W(x,x)<\infty$.
\end{romenumerate}
Note that  these conditions are
satisfied if $W$ is integrable, but they are also satisfied for some
non-integrable $W$. 

A \emph{graphex} is a triple $\cW=(I,S,W)$, where $I\ge0$ is a non-negative
real number, $S:\bbRp\to\bbRp$ is measurable with
$S\land 1$ integrable and $\cW$ is a graphon.
Let $\WW$ be the set of all graphexes.

Each graphex $\cW\in\WW$ defines a random adjacency measure $\xi=\xi(\cW)$, and
thus a corresponding random labelled graph $\gG=\gG(\cW)$, by the
following
construction; 
%(similar to the construction of a random graph sequence from a classical
%graphon on $\oi$); 
see further \citet{VR,VR2} 
and \citet{Kallenberg2}:
Take  realizations of independent unit-rate Poisson processes 
$\Xi=\set{(\gth_j,\vth_j)}_j$ on $\bbRpp$,
$\Xi'_i=\set{(\gs_{ij},\chi_{ij})}_j$ on $\bbRpp$ for $i\in\bbN$,
and $\Xi''=\set{(\rho_j,\rho'_j,\eta_j)}_j$ on $\bbRp^3$.
We regard $\gth_i$, $\gs_{ij}$, $\rho_i$ and $\rho_j'$ as potential vertex
labels, while $\vth_j$, $\chi_{ij}$ and $\eta_j$ can be regarded as types of
the corresponding labels.
Given $\cW=(I,S,W)$ and
these realizations, and a family of \iid{} random variables
$\zeta_{i,j}\sim U(0,1)$ independent of them, define
the adjacency measure
\begin{equation}
  \begin{split}
\xi(\cW)&=\sum_{i,j}\ett{\zeta_{i,j}\le W(\vth_i,\vth_j)  }\gd_{\gth_i,\gth_j}
\\&\qquad
+\sum_{j,k}\ett{\chi_{jk}\le
  S(\vth_j)}\bigpar{\gd_{\gth_j,\gs_{jk}}+\gd_{\gs_{jk},\gth_j}}
\\&\qquad
+\sum_{k}\ett{\eta_k\le I}\bigpar{\gd_{\rho_k,\rho'_k}+\gd_{\rho'_k,\rho_k}}.
  \end{split}
\end{equation}
In other words, the corresponding random labelled graph $\gG(\cW)$ is defined to
have the following edges, with all random choices independent:
\begin{enumerate}
 \renewcommand{\labelenumi}{\textup{(G\arabic{enumi})}}%
 \renewcommand{\theenumi}{\labelenumi}
\item   \label{gG1}
$(\gth_i,\gth_j)$ with probability $W(\vth_i,\vth_j)$ for each pair $(i,j)$
with $i\le j$,
\item \label{gG2}
$(\gth_j,\gs_{jk})$ for each $j$ and $k$ with  $\chi_{jk}\le S(\vth_j)$
\item \label{gG3}
$(\rho_k,\rho_k')$ for each $k$ such that $\eta_k\le I$.
\end{enumerate}
Equivalently, we can define $\gG$ by starting
with the Poisson process $\Xi=(\gth_j,\vth_j)$
and first define the edges in \ref{gG1}, and then add for each $j$ a star with
centre in $\gth_j$ and peripheral vertices labelled by a Poisson process
$\set{\gs_{jk}}_k$ on $\bbRp$ with intensity $S(\vth_j)$,
and finally add edges $(\rho_k,\rho'_k)$ according to a Poisson process with
intensity $2I$ in $\set{(x,x')\in\bbRpp:x<x'}$. (Again all random choices
are independent.)

It follows from more general results by 
\citet{Kallenberg1990}, \cite[Theorem 9.24]{Kallenberg2}
that this construction yields all jointly exchangeable random adjacency
measures, provided we allow the graphex $\cW$ to be random;
see \citet{VR} for the present context.
(See also \cite[Theorem 2.21]{BCCH16} for a related result.)
The conditions \ref{graphona}--\ref{graphono} are precisely the
conditions
needed to guarantee that the constructed measure $\xi$ \as{} is locally
finite (and thus an adjacency measure), see \cite{VR} and
\cite[Proposition 9.25]{Kallenberg2}.

Note that the edges of type \ref{gG3} are independent of everything else and
\as{} isolated; they form a dust of little interest.
Also the stars produced by \ref{gG2} are of minor interest.
One therefore often takes $S=0$ and $I=0$.
A graphex $(0,0,W)$ can be identified with the graphon $W$,
and we write
$\xi(W)=\xi(0,0,W)$ and $\gG(W)=\gG(0,0,W)$. 
In this case, the construction uses only the Poisson process
$(\gth_j,\vth_j)$ and gives the edges in \ref{gG1}; 
see also \cite{BCCH16} (at least when $W$ is integrable).

We write for convenience $\gG_s(\cW)=\gG(\cW)\restr{s}$, and we are particularly
interested in the graph valued process $(\gG_s(\cW))_{s\ge 0}$. Note that this is
an increasing process of finite labelled graphs, where $\gG_0$ is empty and
$\gG_r$ is an induced subgraph
of $\gG_s$ whenever $0\le r\le s$.
Furthermore, $\gG=\bigcup_{s\ge0}\gG_s$, so the (typically infinite) graph
$\gG$ (or the measure $\xi(\gG)$) and the process $(\gG_s)_{s\ge0}$ 
determine each other. 

We consider also the corresponding processes
$\xi_s(\cW):=\xi(\gG_s(\cW))=\xi(\cW)\restr{s}$ of finite adjacency measures and
$G_s(\cW):=\cG(\gG_s(\cW))$ of finite unlabelled graphs.
It follows from the construction above, that for every fixed $s>0$, given the 
unlabelled graph $G_s(\cW)$, the vertex labels on the labelled graph
$\gG_s(\cW)$ are \iid{} with the distribution $U(0,s)$; in other words, conditioned on
$G_s(\cW)$, and therefore also unconditionally,
\begin{equation}\label{lbl=}
  \gG_s(\cW)\eqd \lbl_s(G_s(\cW)).
\end{equation}

\subsection{The classical case of graphons on $\oi$}\label{SSclassical}

As said above, the classical theory of graphons considers graphons
$W:\oii\to\oi$ defined on $\oii$; we can identify them with graphons defined
on $\bbRpp$ that vanish outside $\oii$. As said above, we furthermore
identify the graphon $W$ and the graphex $(0,0,W)$. (There is no dust and no
added stars in the classical theory.) 

In the classical theory, one defines 
for each $n\in\bbN$ 
a random graph $\GG_n(W)$ with $n$ vertices
by taking $n$ \iid{}  random numbers $\vth_i\sim U(0,1)$ and conditionally on
these variables, letting there be an edge $ij$ with probability
$W(\vth_i,\vth_j)$ for each pair $ij$. 
(Furthermore, there are no loops. We may impose this by assuming 
that $W$ vanishes on the diagonal $\gD:=\set{(x,x)}_{x\in\bbRp}$.
For convenience, we tacitly assume this; note that redefining 
a graphon $W$ to be 0 on $\gD$ is equivalent to ignoring loops.)

On the other hand, in the construction above of $\gG_s(W)$, 
\ref{gG2} and \ref{gG3} do not appear, and for \ref{gG1} we can ignore every
$(\gth_i,\vth_i)$ with $\vth_i>1$ (since then $W(\vth_i,\vth_j)=0$ for every
$j$).
In the construction of $\gG_s(W)$ and $G_s(W)$, we thus consider points
$(\gth_i,\vth_i)\in\os\times\oi$. Let the number of these points be
$N_s:=\Xi(\os\times\oi)\sim\Po(s)$. 
It follows that conditioned on $N_s$, the random graph $G_s(W)$ constructed
above equals $\GG_{N_s}(W)$, with all isolated vertices deleted.

\section{Graphex equivalence and convergence}\label{Sconv}

In this section, we define the four different types of convergence that we
consider in the present paper.

\subsection{Convergence of  corresponding random graphs: $\togp$}
\label{Stogp}
The construction in \refS{Snot} defines a random labelled graph $\gG(\cW)$ and a
corresponding random adjacency measure $\xi(\cW)$ for every graphex $\cW$.
This defines a map $\Psi:\WW\to\PMRR$ by $\Psi(\cW):=\cL(\xi(\cW))$, the
distribution of $\xi(\cW)$. Unfortunately, $\Psi$ is not injective, \ie, a
graphex $\cW$  is not uniquely determined by the distribution of $\xi(\cW)$.
We say that two graphexes $\cW$ and $\cW'$ are \emph{equivalent} if
$\cL(\xi(\cW))=\cL(\xi(\cW'))$; we denote this by $\cW\cong\cW'$.
Let $\WWX:=\WW/\cong$, the set of equivalence classes; 
then $\Psi$ can be regarded as an injection $\WWX\to\PMRR$, and we can
identify $\WWX$ with its image and equip $\WWX$ with the subspace topology
inherited from $\PMRR$.
Since $\MRR$ is a Polish space,
$\PMRR$ is metrizable, and thus the topology of $\WWX$ can be defined by a
metric $\dgp$. (There are many possible choices $\dgp$ 
but we assume that one is chosen; we will not distinguish an
explicit choice.) 
Moreover, we can also regard $\dgp$ as defined on $\WW$; this makes $\WW$
into a pseudometric space, with 
\begin{equation}\label{equiv}
  \cW\cong\cW' \iff \dgp(\cW,\cW')=0.
\end{equation}

Convergence in the pseudometric $\dgp$ is denoted $\togp$ by
\citet{VR2}. (They actually use \ref{T0G} in \refT{T0} below as the definition.)
Convergence $\togp$
can be characterized as follows, which is at least implicit in \cite{VR2} 
but stated explicitly here for easy reference.
Note that since the topology is metric, it suffices to consider convergence
of sequences; the theorem extends immediately to, \eg, convergence of
families with a continuous parameter.

\begin{theorem}\label{T0}
  Let $\cW_n$, $n\ge1$, and $\cW$ be graphexes.
Then the following are equivalent,  as \ntoo{}. 
\begin{romenumerate}
\item \label{T0gp}
$\cW_n\togp\cW$.
\item \label{T0dgp}
$\dgp(\cW_n,\cW)\to0$.
\item \label{T0cL}
$\cL(\xi(\cW_n)) \to \cL(\xi(\cW))$ in $\PMRR$.
\item \label{T0dxi}
$\xi(\cW_n)\dto\xi(\cW)$ in $\MRR$.
\item \label{T0xis}
$\xi_s(\cW_n) \dto\xi_s(\cW)$ in $\MRR$
for every $s<\infty$.
\item \label{T0G}
$G_s(\cW_n) \dto G_s(\cW)$ in $\fGf$, for every $s<\infty$.
\end{romenumerate}
In \ref{T0xis}--\ref{T0G}, it suffices to consider $s$ in a given
unbounded subset of $\bbRp$, for example $s\in\bbN$.
\end{theorem}

\begin{proof}
\ref{T0gp}$\iff$\ref{T0dgp}$\iff$\ref{T0cL} holds by the definitions above.

\ref{T0cL}$\iff$\ref{T0dxi} holds by the definition of convergence in
distribution. 

\ref{T0dxi}$\implies$\ref{T0xis}.
Let,
$\chS$ be the class of bounded measurable subsets of $\bbRpp$, and for a
random measure $\xi$, let 
$\chS_\xi:=\set{B\in\chS:\xi(\partial B)=0 \text{ \as}}$.
By \cite[Theorem 16.16]{Kallenberg}, a sequence of random measures
$\xi_n\dto\xi$ if and only if, for every finite sequence
$B_1,\dots,B_k\in\chS$,
\begin{equation}\label{b16}
  \bigpar{\xi_n(B_1),\dots,\xi_n(B_k)}
\dto   \bigpar{\xi(B_1),\dots,\xi(B_k)}.
\end{equation}

Let $Q_s:=[0,s]^2$. The random graph $\gG(\cW)$ has \as{} no label $s$, and thus
the random measure $\xi(\cW)$ has \as{} no mass at $\partial Q_s$, \ie,
$Q_s\in\chS_{\xi(\cW)}$.
Hence, if 
$B\in\chS_{\xi_s(\cW)}$, then
$B\cap Q_s\in\chS_{\xi(\cW)}$.
It follows, see \eqref{b16}, 
that if $\xi(\cW_n)\dto\xi(\cW)$ and $B_1,\dots,B_k\in\chS_{\xi_s(\cW)}$, then
\begin{equation}
  \begin{split}
    \bigpar{\xi_s(\cW_n)(B_j)}_{j=1}^k
&=
    \bigpar{\xi(\cW_n)(B_j\cap Q_s)}_{j=1}^k
\\&
\dto
    \bigpar{\xi(\cW)(B_j\cap Q_s)}_{j=1}^k
=
    \bigpar{\xi_s(\cW)(B_j)}_{j=1}^k. 
  \end{split}
\end{equation}
Consequently, $\xi_s(\cW_n)\dto\xi_s(\cW)$ by  
\cite[Theorem 16.16]{Kallenberg} again.
(Alternatively, this implication follows by the continuous mapping theorem since
a similar argument shows that the mapping $\xi\to\xi\restr s$ from $\MRR$ to itself is continuous
at every $\xi$ with $\xi(\partial Q_s)=0$.)

\ref{T0xis}$\implies$\ref{T0dxi}. 
Let $B_1,\dots,B_k\in\chS_{\xi(\cW)}$, and let $s$ be so large that  $\bigcup_j
B_j\subset Q_s$. Then $\xi_s(B_j)=\xi(B_j)$ and thus
$\xi_s(\cW_n)\dto \xi_s(\cW)$ implies, 
again using \cite[Theorem 16.16]{Kallenberg},
$\bigpar{\xi(\cW_n)(B_j)}_{j=1}^k\dto    \bigpar{\xi(\cW)(B_j)}_{j=1}^k$.
Hence, $\xi(\cW_n)\dto\xi(\cW)$.

\ref{T0xis}$\iff$\ref{T0G}.
Fix $s>0$.
By \eqref{lbl=}, $\xi_s(\cW)\eqd\xi(\lbls(G_s(\cW)))$, and similarly for
$\cW_n$. Hence the equivalence follows from 
\refL{LA} below, taken from \cite{VR2}.
%\cite[Lemma 4.11]{VR2}. 
\end{proof}

\begin{lemma}[{\cite[Lemma 4.11]{VR2}}]\label{LA}
If $G_n$, $n\ge1$, and $G$ are  any random graphs in   $\fGf$, and $s>0$,
then $G_n\dto G$ in $\fGf$
if and only if $\lbls(G_n)\dto\lbls(G)$ in $\fGl$, \ie,
if and only if $\xi(\lbls(G_n))\dto\xi(\lbls(G))$ in $\MRR$.
\end{lemma}

We will not repeat the proof of \cite[Lemma 4.11]{VR2}, but we note that it
can be interpreted as defining a map 
$\Phi_s:\cP(\fGf)\to\cP(\MRR)$, by taking for $\mu\in\cP(\fGf)$ a random
graph $G\sim\mu$, and defining
$\Phi_s(\mu):=\cL(\xi(\lbl_s(G)))\in\cP(\MRR)$.
It is then shown that $\Phi_s$ is continuous, injective and proper.
Finally, any continuous and proper map to a metric space is closed, and a
continuous and closed injection is a homeomorphism onto a closed subset.

\subsection{Cut metric: $\dcut$}

The invariant cut metric $\dcut$ is defined only for integrable graphons.
It is the standard metric for classical graphons, see \eg{} \cite{BCLSV1},
\cite{Lovasz}, \cite{SJ249}. The definition was extended to
graphons defined on arbitrary $\gs$-finite
measure spaces by \citet{BCCH16}, to which we refer for details.
Here we only consider graphons on $\bbRp$, and then the results simplify as
follows, see \cite{BCCH16}.

First, for an integrable function $F$ on $\bbRpp$, we define its 
\emph{cut  norm} by
\begin{equation}\label{cn}
  \cn{F}
:= \sup_{T,U}\lrabs{\int_{T\times U} F(x,y)\dd\mu(x)\dd\mu(y)},
\end{equation}

If $W$ is a graphon and $\gf:\bbR\to\bbR$ is measure preserving,
let $W^\gf(x,y):=W(\gf(x),\gf(y))$.
Then, 
for two integrable graphons $W_1$ and $W_2$,  define
\begin{equation}\label{dcut1}
  \dcut(W_1,W_2) :=\inf_{\gf_1,\gf_2}\cn{W_1^{\gf_1}-W_2^{\gf_2}},
\end{equation}
taking the infimum over all pairs of measure preserving maps
$\gf_1,\gf_2:\bbRp\to\bbRp$.
Moreover, it is shown in \cite[Proposition 4.3(c)]{BCCH16}
that
\begin{equation}\label{dcutR}
  \dcut(W_1,W_2) :=\inf_{\gf}\cn{W_1^\gf-W_2},
\end{equation}
taking the infimum over all \mpb{s}  $\gf:\bbRp\to\bbRp$.

By \cite[Theorem 2.22]{BCCH16}, $\dcut(W_1,W_2)=0$ if and only if $W_1\cong
W_2$, \ie, $W_1$ and $W_2$ are equivalent in the sense defined above.
(For characterizations of this, see \cite{BCCH16} and \cite{SJ311}.)

If $W_n$ and $W$ are integrable graphons, we write $W_n\tocut W$ as \ntoo{}
if
$\dcut(W_n,W)\to 0$.
We shall show that this is at least as strong as $\togp$.
We begin with the special case of classical graphons, where we have equivalence.

\begin{lemma}\label{LF}
Suppose that $W_n$, $n\ge1$, and $W$  are graphons with support on $\oii$.
Also suppose that they all vanish on the diagonal $\gD:=\set{(x,x)}_{x\in\oi}$.
Then the following are equivalent, as \ntoo:
\begin{romenumerate}
\item \label{LFcut}
$\dcut(W_n,W)\to 0 $

%\item for some $s>0$,
%$G_s(W_n)\dto G_s(W)$.

\item \label{LFs}
for every $s\ge0$,
$G_s(W_n)\dto G_s(W)$.

\item \label{LFgp}
$W_n\togp W$.
\end{romenumerate}
\end{lemma}

\begin{proof}
\ref{LFcut}$\implies$\ref{LFs}.
Using the notation in \refS{SSclassical},  \ref{LFcut} implies (and is
equivalent to)
$\GG_N(W_n)\dto \GG_N(W)$ as \ntoo{} for every $N\in\bbN$, see \eg{}
\cite{BCLSV1,SJ209}. Hence the same holds if we let $N=N_s\sim\Po(s)$ be
random, and since the space of finite graphs is discrete, the 
result holds also if we remove all isolated vertices from the graphs, which
yields \ref{LFs}, see \refS{SSclassical}.

\ref{LFs}$\iff$\ref{LFgp}.  By \refT{T0}.

\ref{LFs}$\implies$\ref{LFcut}. 
The implication \ref{LFcut}$\implies$\ref{LFs} shows 
that the map $W\mapsto (\cL(G_s(W)))_{s\ge0}\in \cP(\fGf)^{\bbRp}$ is
continuous.
Since the space of (equivalence classes) of classical graphons is compact,
it suffices to show that the map is injective.

The only complication is caused by the removal of isolated vertices.
Thus, let $H$ be a finite graph without isolated vertices, and let $H+mK_1$
be $H$ with $m$ isolated vertices added.
Then, with $v:=v(H)$,
\begin{equation}\label{pow}
  \begin{split}
\P(G_s(W)=H) 
&= \summo\P\bigpar{\GG_{N_s}(W)=H+mK_1}
\\&
=e^{-s} \sum_{m=0}^\infty \frac{s^{m+v}}{(m+v)!}\P\bigpar{\GG_{m+v}(W)=H+mK_1}.
\end{split}
\end{equation}
Hence, if $W$ and $W'$ are two classical graphons and
$\P(G_s(W)=H)=\P(G_s(W')=H)$, then multiplication by $e^s$ and
identification of the coefficients of the power series in \eqref{pow}
show that $\P(\GG_{m+v}(W)=H+mK_1)=\P(\GG_{m+v}(W')=H+mK_1)$ for every
$m\ge0$.
Consequently, $G_s(W)\eqd G_s(W')$ for every $s>0$ implies that
$\GG_N(W)\eqd \GG_N(W')$ for every $N$, and thus $W\cong W'$ and
$\dcut(W,W')=0$. 
\end{proof}

\begin{theorem}\label{T4}
  Suppose that $W_n$, $n\ge1$, and $W$  are integrable graphons.
If $W_n\tocut W$ as \ntoo, then $W_n\togp W$.
\end{theorem}

\begin{proof}
 Assume $\dcut(W_n,W)\to0$.
By \eqref{dcutR}, we may then replace each $W_n$ by an equivalent
$W_n^{\gf_n}$
such that $\cn{W_n-W}\to0$; note that the measurable rearrangement
$W_n^{\gf_n}$ defines the same random graphs as $W_n$ in distribution, \ie,
$\xi(W_n^{\gf_n})\eqd \xi(W_n)$, and thus $\dgp(W_n,W_n^{\gf_n})=0$, see
\eqref{equiv}. 

In the sequel we thus assume $\eps_n:=\cn{W_n-W}\to0$.
Define the truncations $W_n\NN:=W_n\etta_{[0,N]^2}$ and $W\NN:=W\etta_{[0,N]^2}$.
Let $\eps>0$ and
fix a large $N$ such that $\norm{W\NN-W}_{L^1}<\eps$.
Then, 
\begin{equation}
  \cn{W_n\NN-W\NN}\le \cn{W_n-W}=\eps_n\to0
\end{equation}
as \ntoo. Consequently, by \refL{LF} and a rescaling,
for every $s\ge0$,
\begin{equation}\label{qqq}
G_s(W_n\NN)\dto G_s(W\NN)
\qquad\text{as \ntoo}.
\end{equation}

Furthermore, 
\begin{equation}\label{kkk}
  \begin{split}
  \int(W_n-W_n\NN)\dd\gl  
&=
\int(W-W\NN)\dd\gl + 
\int_{\bbRpp\setminus[0,N]^2}(W_n-W)\dd\gl
\\&
<\eps + 2 \cn{W_n-W} =\eps+2\eps_n.    
  \end{split}
\end{equation}

$G_s(W)$ and $G_s(W\NN)$ differ only if the 
labelled graph $\gG_s(W)$ contains some edge with at least one label $>N$.
The expected number of such edges is
\begin{equation}
\E \bigpar{e(G_s(W))-e(G_s(W\NN))}
=\frac{s^2}2 \int_{\bbRpp} (W-W\NN)\dd\gl,
\end{equation}
 and thus
\begin{equation}
  \P\bigpar{G_s(W)\neq G_s(W\NN)} \le 
s^2 \int_{\bbRpp} (W-W\NN)\dd\gl < s^2\eps,
\end{equation}
and similarly, using \eqref{kkk},
\begin{equation}\label{ccc}
  \P\bigpar{G_s(W_n)\neq G_s(W_n\NN)} \le 
 s^2\int_{\bbRpp} (W_n-W_n\NN)\dd\gl < s^2(\eps+2\eps_n).
\end{equation}

Consequently, if $f:\fGf\to\oi$ is any function, then by \eqref{qqq} and
\eqref{ccc}, 
\begin{equation*}
  \begin{split}
  \limsup_\ntoo&|\E f(G_s(W_n))-\E f(G_s(W))|
\\&
\le
  \limsup_\ntoo\bigpar{|\E f(G_s(W_n\NN))-\E f(G_s(W\NN))| + s^2(2\eps+2\eps_n)}
\\&
= 2s^2\eps.    
  \end{split}
\end{equation*}
Since $\eps>0$ is arbitrary, this shows 
$\E f(G_s(W_n))\to\E f(G_s(W))$, and thus $G_s(W_n)\dto G_s(W)$ as \ntoo.
The result follows by \refT{T0}.
\end{proof}

\begin{remark}
In \refL{LF}\ref{LFs}, it is not necessary to assume the condition 
$G_s(W_n)\dto G_s(W)$ for every
$s\ge0$. In fact, 
it suffices to assume this for $s$ in an arbitrary  non-empty
interval $(a,b)$, or even more generally, for $s$ in any infinite set having
a cluster point in $\ooo$; this follows by the same proof and the uniqueness
theorem for analytic functions.

Moreover, the argument suggests that it might suffice to assume \ref{LFs}
for a single $s>0$. We state this as an open problem. 
By the proof above, this is equivalent to the problem
whether $G_s(W)\eqd G_s(W')$ for the random graphs without isolated vertices
is equivalent to the corresponding equality in distribution of the random
graphs $\GG_{\Po(s)}(W)$ and $\GG_{\Po(s)}(W')$ with
isolated vertices.
\end{remark}

\begin{problem}
In \refL{LF}, does $G_s(W_n)\dto G_s(W)$ for a single $s>0$ imply $W_n\togp W$? 
\end{problem}

We can also ask whether this holds for general graphexes and not
just for classical graphons. We leave this too as open problems.

\begin{problem}
In \refT{T0}, 
does $G_s(\cW_n)\dto G_s(\cW)$ for, say, $0<s<1$ imply $\cW_n\togp \cW$?   
\end{problem}

\begin{problem}
In \refT{T0}, 
does $G_s(\cW_n)\dto G_s(\cW)$ for a single $s>0$ imply $\cW_n\togp \cW$?   
\end{problem}

\subsection{Stretched convergence: $\dcuts$ and $\togs$}
\label{SSstretch}

We define, as in \cite{BCCH16} and \cite{VR2}, given a graphon $W$ or more
generally a graphex $\cW=(I,S,W)$ and a real number $c>0$, the
\emph{stretched} graphon or graphex by
$W\strc(x,y):=W(x/c,y/c)$ or $\cW\strc:=(c^2I,S\strc,W\strc)$ where further
$S\strc(x):=cS(x/c)$.  It follows easily from the construction of the random
graphs above that a stretched graphex defines the same random graph process
$G_r(\cW)$ up to a change of parameter:
\begin{equation}\label{para}
(G_r(\cW\strc))_r\eqd (G_{cr}(\cW))_r.  
\end{equation}

\citet{BCCH16} define the \emph{stretched cut metric} $\dcuts$
by, for two nonzero integrable graphons $W_1$ and $W_2$,  
\begin{equation}
  \dcuts(W_1,W_2):=\dcut\bigpar{W_1\str{c_1},W_2\str{c_2}}
\end{equation}
where $c_i:=\norm{W_i}_{L_1}\qqw$. 
We write $W_n\tocuts W$ as \ntoo{} if
$\dcuts(W_n,W)\to 0$.
Since $\dcut(W_n,W)\to0$ implies
$\norm{W_n}_{L^1}\to\norm{W}_{L^1}$, it is easily seen that, for any
integrable graphons $W_n$ and $W$,
\begin{equation}\label{dcuts}
  W_n\tocuts W \iff W_n\str{c_n}\tocut W\strc 
\text{ for some constants $c_n,c>0$}.
\end{equation}

Moreover, note that the random graph process 
$r\mapsto \gG_r(\cW)$ of labelled graphs is
increasing and right-continuous (and thus cadlag), and has a.s.\ only
a finite number of jumps $\tau_k$ in each finite interval.
Hence the unlabelled graph process
$r\mapsto G_r(\cW)=\cG(\gG_r(\cW))$ 
has the same properties.
\citet{VR2} consider for any non-zero graphex
$\cW$ the sequence $(G_{\tau_k}(\cW))_k$ of different 
(finite and unlabelled) graphs that appear in $\set{G_r(\cW):0\le r<\infty}$; 
they define
for $\cW_n,\cW\in\WW':=\WW\setminus\set0$,
\begin{equation}
  \cW_n\togs\cW \qquad \text{if } 
(G_{\tau_k}(\cW_n))_{k=1}^\infty \dto (G_{\tau_k}(\cW))_{k=1}^\infty.
\end{equation}

This too is a metric convergence.
In analogy with the  definition of $\dgp$ above, we can define a map
$\Psi':\WW\setminus\set0\to \cP(\fGf^\infty)$  by
$\Psi(\cW)=\cL\bigpar{(G_{\tau_k}(\cW))_{k\ge1}}$, and by fixing a metric on  
$\cP(\fGf^\infty)$, we define a pseudometric $\dgs$ on
$\WW'$.
Then obviously
\begin{equation}
  \cW_n\togs \cW \iff \dgs(\cW_n,\cW)\to0.
\end{equation}

It follows by \eqref{para} that $(G_{\tau_k}(\cW))_{k}\eqd
(G_{\tau_k}(\cW\strc))_{k}$ for any stretching of a graphex $\cW$, and thus 
$\dgs(\cW,\cW\strc)=0$, \cf{} \cite[Corollary 5.5]{VR2}.
Furthermore, it follows from \refT{T0} that, assuming $\cW,\cW\in\WW'$,
see \cite[Lemma 5.6]{VR2} for a detailed proof,
\begin{equation}
  \cW_n\togp \cW \implies \cW_n\togs\cW.
\end{equation}

Consequently, we have the following partial analogue of \eqref{dcuts}.
\begin{lemma}[\cite{VR2}]\label{Ltogs}
If $\cW_n,\cW\in\WW$ and there exist $c_n,c>0$ such that 
$\cW_n\str{c_n} \togp \cW\strc$, then 
$\cW_n\togs \cW$.
\qed
\end{lemma}

\begin{problem}
  Does the converse to \refL{Ltogs} hold?
(We conjecture so.)
\end{problem}

We also have a result corresponding to \refT{T4} for the stretched metrics.

\begin{theorem}\label{T5}
  Suppose that $W_n$, $n\ge1$, and $W$  are integrable non-zero graphons.
If $W_n\tocuts W$ as \ntoo, then $W_n\togs W$.
\end{theorem}
\begin{proof}
  By \eqref{dcuts}, we have 
$W_n\str{c_n} \tocut W\strc$ for some constants $c_n$ and $c$, and the
result follows by \refT{T4} and \refL{Ltogs}.
\end{proof}

\begin{problem}
Does the converses of Theorems \ref{T4} and \ref{T5} hold.
(We conjecture not.)
\end{problem}

\section{Random relabellings}

We show in this section that \cite[Theorem 4.3]{VR2} extends 
to convergence as \stoo{} through the set of all positive real numbers.
%We state this result as follows.

\begin{theorem}[{Extension of \cite[Theorem 4.3]{VR2}}] \label{T1}
Let $\cW$ be a graphex and let $G_s:=\cG(\gG_s(\cW))$, $s\ge0$, 
be the corresponding process of unlabelled graphs.  
Then, \as, as \stoo,
$\bigpar{\lbl_s(G_s)\mid G_s}\dto \gG(\cW)$ in $\fGl$, \ie,
$\cL(\lbl_s(G_s)\mid G_s)\to \cL(\gG(\cW))$ in $\cP(\fGl)$.
Equivalently,
$\bigpar{\xi(\lbl_s(G_s))\mid G_s}\dto \xi(\cW)$ in $\MRR$.
%$\cL(\lbl_s(G_s)\mid G_s)\to \cL(\gG(\cW))$ in $\PMRR$, 
\end{theorem}

Note that \refT{T1} not only generalizes \cite[Theorem 4.3]{VR2}; 
also \cite[Theorem 4.4]{VR2} is an immediate corollary.

In the proof we  use the following lemma, which extends a standard lemma
used by \cite{VR2} from a parameter $s\in \bbN$ to a parameter $s\in\bbQp$. 
We guess that also the version here is known, but since we do not know a
reference, 
we give a simple proof for completeness.

\begin{lemma}\label{LT1}
Assume that $X_s$ is a random variable for each $s\in\bbQ_+$ such that
$X_s\to X$ \as{} as \stoo, and further \as{} $|X_s|\le Y$ for some random
variable $Y$ with $\E Y<\infty$.
Let $(\cF_s)_{s\in\bbQp}$  be a decreasing family of $\gs$-fields and let
$\cF_\infty:=\bigcap_{s\in\bbQp}\cF_s$. 
Then $\E\bigpar{X_s\mid\cF_s}\to\E\bigpar{X\mid\cF_\infty}$ \as{} as
  \stoo{} with $s\in\bbQp$.
\end{lemma}

\begin{proof}
%Let $Z_t:=\sup_{s\ge t, s\in\bbQp}|X_s-X|$,
Let $Z_t:=\sup_{s\ge t}|X_s-X|$,
where as in the rest of the proof
we consider only $s\in \bbQp$.
The assumption implies \as{} $|X|\le Y$ and
thus $0\le Z_t\le 2Y$.
Moreover, $Z_t\to0$ \as{} as \ttoo{} by assumption. 
Hence, $\E Z_t<\infty$ and $\E Z_t\to0$ as \ttoo{} by
dominated convergence.

Fix $t$. For $s\ge t$ we have $|X_s-X|\le Z_t$, and thus
$\E\bigpar{|X_s-X|\mid\cF_s}\le \E\bigpar{Z_t\mid\cF_s}$ a.s.
Consequently, using the convergence theorem for 
(reverse) martingales, 
\cite[Theorem 7.23]{Kallenberg},
\as,
\begin{equation}
\limsup_{s\to\infty} \E\bigpar{|X_s-X|\mid\cF_s}
\le\limsup_{s\to\infty} \E\bigpar{Z_t\mid\cF_s}
=\E\bigpar{Z_t\mid\cF_\infty}.
\end{equation}
Hence, for every $t\ge0$,
\begin{equation}\label{ragna}
\E\limsup_{s\to\infty} \E\bigpar{|X_s-X|\mid\cF_s}
\le
\E\bigpar{\E\bigpar{Z_t\mid\cF_\infty}}
=\E Z_t.
\end{equation}
However, we have shown that $\E Z_t\to0$ as \ttoo, and thus \eqref{ragna}
implies 
\begin{equation}%\label{ragna}
\E\limsup_{s\to\infty} \E\bigpar{|X_s-X|\mid\cF_s}=0.
\end{equation}
Consequently, \as, 
$\E\bigpar{|X_s-X|\mid\cF_s}\to0$ as \stoo, and thus
\begin{equation}\label{eva}
\bigabs{\E\bigpar{X_s-X\mid\cF_s}}\le
\E\bigpar{|X_s-X|\mid\cF_s}\to0.  
\end{equation}

Furthermore, by the (reverse) martingale convergence theorem again and
\eqref{eva}, a.s.,
\begin{equation}
  \E\bigpar{X_s\mid\cF_s}
= \E\bigpar{X\mid\cF_s} + \E\bigpar{X_s-X\mid\cF_s}
\to \E\bigpar{X\mid\cF_\infty} + 0 .
\end{equation}
\end{proof}

\begin{proof}[Proof of \refT{T1}]
  Let $(\tau_k)_k$ be the jump times, where $G_s\neq G_{s-}$.
Then $G_s$ is constant for $s\in[\tau_k,\tau_{k+1})$, and it is easily seen
that the distribution
$\cL(\lbl_s(G_s)\mid G_s)$ is a continuous function of $s\in[\tau_k,\tau_{k+1})$.
It follows that it suffices to prove the result as $s\to\infty$ through the
countable set of rational numbers.
Thus, in the sequel of the proof, we assume that $s\in\bbQ_+$, and we
consider limits as $\stoo$ through $\bbQp$.

Except for this, we follow the proof of \cite[Theorem 4.3]{VR2}; 
for completeness we repeat most of the arguments.

First, by \cite[Lemma 4.2]{VR2}, a consequence of 
\cite[Theorems 16.28--29]{Kallenberg},
it suffices to prove that if
$U$ is a finite union of rectangles with rational coordinates, 
then, with $\xi_s:=\xi(\lbl_s(G_s))$ and $\xi:=\xi(\cW)$,
\begin{equation}\label{ga}
\cL(\xi_s(U)\mid G_s)\to \cL\bigpar{\xi(U)}  
\qquad\text{a.s.}
\end{equation}
as \stoo.
Note that \cite[Lemma 4.2]{VR2}, although stated for sequences, immediately
extends to a parameter $s\in\bbQp$ (or even $\bbRp$) since convergence in
$\PMRR$ can be defined by a metric.

Next, \eqref{ga} means that for every
function $f:\bbZp\to \bbRp$ of the form 
$f(x)=\etta_{\set{k}}(x)=\ett{x=k}$,
\begin{equation}\label{gb}
 \E\bigpar{f(\xi_s(U))\mid G_s}\to\E f(\xi(U))
\qquad \text{a.s.}
\end{equation}

Thus, fix such a $U$ and $f$ (or any bounded $f:\bbZp\to\bbRp$),
and fix some large enough $r$ such that $U\subset [0,r]^2$.
For $s\in\bbQp$, let $\gG^s$ be the partially labelled graph obtained from
$\gG(\cW)$ by forgetting all labels in $\os$ (but keeping larger labels).
Let $\cF_s$ be the $\gs$-field generated by $\gG^s$.
Conditioned on $\gG^s$, $\gG(\cW)$ is obtained by randomly relabelling the
unlabelled vertices by labels in $\os$; 
note that the unlabelled part of $\gG^s$ is $G_s$.
Hence, if $s>r$, 
\begin{equation}\label{eff}
  \E\bigpar{ f(\xi(U))\mid \cF_s}=
\E\bigpar{f(\xi(\lbls(G_s))\mid G_s}=\E\bigpar{f(\xi_s(U))\mid G_s}.
\end{equation}

Define
as in \cite{VR2} $U_t:=U+(t,t)$ and, for $s>r$,
\begin{equation}\label{xrs}
  X_s\rr:=\frac{1}{s-r}\int_0^{s-r} f(\xi(U_t))\dd t.
\end{equation}
Since $t\mapsto Y_t:=f(\xi(U_t))$ is a stationary stochastic process, which
is $r$-dependent in the sense that $\set{Y_t}_{t\le t_0}$ is independent of
$\set{Y_t}_{t>t_0+r}$ for every $t_0$, it follows from the ergodic theorem 
\cite[Corollary 10.9]{Kallenberg} that $X_s\rr\asto \E f(\xi(U))$ as \stoo.
\refL{LT1} now shows that 
\begin{equation}\label{xr1}
\E\bigpar{ X_s\rr\mid\cF_s}\asto \E f(\xi(U)).
\end{equation}
Furthermore, by symmetry, for every $s$ and every $t<s-r$,
$\E\bigpar{f(\xi(U_t))\mid \cF_s} = \E\bigpar{f(\xi(U))\mid \cF_s}$ \as,
and thus \eqref{xrs} implies
\begin{equation}\label{xr2}
\E\bigpar{X_s\rr\mid \cF_s} = \E\bigpar{f(\xi(U))\mid \cF_s} 
\qquad \text{a.s.}
\end{equation}

We obtain \eqref{gb} by combining \eqref{eff}, \eqref{xr2} and \eqref{xr1},
which completes the proof.
\end{proof}

\section{Convergence of empirical graphons}

If $G$ is any finite graph, and $s>0$, we define as in \cite{VR2} the
\emph{stretched empirical graphon} $\hW_{G,s}$ as 
the graphon obtained from 
the adjacency matrix of $G$ by replacing each vertex by an interval of
length $1/s$. (This assumes some ordering of the vertices, but different
orderings give equivalent graphons.) 
In other words, using the notation of \refS{SSstretch},
$\hW_{G,s}=\hW_G\str{1/s}$ where $\hW_G:=\hW_{G,1}$.
Note that $\hW_{G,s}$ only takes the
values \setoi.

We identify as usual the graphon $\hW_{G,s}$ with the graphex $(0,0,\hW_{G,s})$.

\citet[Theorem 2.23]{BCCH16} show that if $W$ is a non-zero integrable
graphon, then 
\as{} $\dcuts\bigpar{\hW_{G_s(W)}, W}\to0$ and (implicitly, or as a consequence)
$\dcut\bigpar{\hW_{G_s(W),s}, W}\to0$ as \stoo.

Similarly, \citet[Theorems 4.8 and 4.12]{VR2} show convergence in $\togp$
for a general graphex; however, only for sequences $s_k\to\infty$.
We show in this section an extension of %\cite[Theorems 4.8 and 4.12]{VR2} 
their result
to convergence as \stoo{} through the set of all positive real numbers.

\begin{theorem}[{Extension of \cite[Theorems 4.8 and 4.12]{VR2}}] \label{T2}
Let\/ $\cW$ be any graphex and let $G_s:=\cG(\gG_s(\cW))$, $s\ge0$, 
be the corresponding process of unlabelled graphs.  
Then, %\as, as \stoo,
$\hW_{G_s,s}\togp \cW$ \as{} as \stoo.
\end{theorem}

\begin{proof}
Again we follow  the proof in \cite{VR2} with some modifications.
Let $0<r\le s$.

Given $G_s$ we define two random induced subgraphs $X\rs$ and $M\rs$ of
$G_s$;
in $X\rs$ we select each vertex of $G_s$ with probability $r/s$, and in
$M\rs$ we take $\Po(r/s)$ copies of each vertex in $G_s$ (thus allowing
repetitions); in both cases this is done independently for all vertices. 
In both cases, we then add the edges induced by $G_s$ and remove any 
isolated vertex.

It follows that, conditionally given $G_s$,
$\lblr(X\rs) \eqd \lbls(G_s)\restrr$.
Hence,
by Theorems \ref{T1} and \ref{T0} together with \eqref{lbl=}, \as{} as \stoo,
\begin{equation*}%\label{erika0}
  \bigpar{\lblr(X\rs)\mid G_s}
\eqd
  \bigpar{\lbls(G_s)\restrr\mid G_s} \dto \gG(\cW)\restrr=\gG_r(\cW)
=\lblr(G_r(\cW)).
\end{equation*}
Consequently, using \refL{LA},
\begin{equation}\label{erika}
  \bigpar{X\rs\mid G_s}
\dto G_r(\cW).
\end{equation}

On the other hand, still conditionally given $G_s$, 
the random graph $G_r(\hWss)$ constructed as in \refS{SSgraphon} from the
empirical graphon $\hWss$ is precisely $M\rs$.
Hence, 
\begin{equation}\label{per}
\bigpar{M\rs\mid G_s}
\eqd \bigpar{G_r(\hWss)\mid G_s}.
%\eqd \bigpar{\gG_r(\hWss)\mid G_s}.
\end{equation}

\refL{LXM} below shows that, still conditionally given $G_s$,
 there exists a coupling of $X\rs$ and $M\rs$ such that
$\P\bigpar{M\rs\neq X\rs\mid X\rs}\le 2(r/s)v(X\rs)$. Since trivially also
the probability is at most 1, this yields, for any constant $A>0$,
$\P\bigpar{M\rs\neq X\rs\mid X\rs}\le 2Ar/s+ \ett{v(X\rs)>A}$
and thus
\begin{equation}\label{ahman}
  \P\bigpar{M\rs\neq X\rs\mid G_s}
%\le \E\Bigpar{1\land 2\frac{r}{s}v(X\rs)\mid G_s}
\le 2Ar/s + \P\bigpar{v(X\rs)>A\mid G_s}.
\end{equation}
Moreover, \eqref{erika} implies that
\begin{equation}\label{mattsson}
  \bigpar{v(X\rs)\mid G_s}
%=  \bigpar{v(\lblr(X\rs))\mid G_s}
\dto v(G_r(\cW))
=
v(\gG_r(\cW))
\end{equation}
\as{} as \stoo,
and thus
\begin{equation}\label{matts}
  \P\bigpar{v(X\rs)>A\mid G_s}
\asto \P\bigpar{v(\gG_r(\cW))>A}.
\end{equation}
It follows from \eqref{ahman} and \eqref{matts}, that for every fixed $r$
and $A$, \as, 
\begin{equation}%\label{ahman}
\limsup_{\stoo}  \P\bigpar{M\rs\neq X\rs\mid G_s}
\le 0 +  \P\bigpar{v(\gG_r(\cW))>A}.
\end{equation}
The \lhs{} does not depend on $A$, and as $A\to\infty$, the \rhs{}
tends to 0. Hence, 
\begin{equation}\label{perp}
  \P\bigpar{M\rs\neq X\rs\mid G_s}\to0
\end{equation}
\as{} as \stoo, for every fixed $r$.
Combining \eqref{perp} with \eqref{erika} and \eqref{per}, we obtain
\begin{equation}
\bigpar{G_r(\hWss)\mid G_s} \dto G_r(\cW)  
\end{equation}
\as{} as \stoo, for every fixed $r$, and thus  \as{} for every $r\in\bbN$.
Consequently, the result follows by \refT{T0}.
% by \refT{T0}, $\hWss\togp\cW$ \as{} as \stoo.
\end{proof}

\begin{lemma}\label{LXM}
Given a finite graph $G$ and $p\in\oi$, let $X_p$ and $M_p$  be the random
induced subgraphs of $G$ obtained by independently taking $\Be(p)$ and
$\Po(p)$ copies of each vertex of $G$, respectively,
and then eliminating all resulting isolated vertices.
Then $X_p$ and $M_p$ may be coupled such that 
$\P(M_p\neq X_p\mid X_p)\le 2p v(X_p)$.
\end{lemma}
\begin{proof}
Let, for $i\in V(G)$, $I_i\sim\Be(p)$ and $Y_i\sim\Po(p)$ be 
random variables with the pairs $(I_i,Y_i)_{i\in V(G)}$ independent, let
$\bX_p$ and $\bM_p$ be the induced subgraphs of $G$ obtained by taking $I_i$ or
$Y_i$ copies of each vertex $i$, respectively, and
let $X_p$ and $M_p$ be the subgraphs of $\bX_p$ and $\bM_p$ obtained by
deleting all isolated vertices.

We couple $I_i$ and $Y_i$ such that for each $i$, $I_i=0\implies Y_i=0$. 
This is possible since $\P(Y_i=0)=e^{-p}\ge1-p=\P(I_i)=0$, 
and we have by a simple
calculation 
\begin{equation}\label{nora}
\P\bigpar{Y_i\neq I_i\mid I_i=1}
=\frac{\P(Y_i\neq I_i)}{\P(I_i=1)}
=
\frac{2(p-pe^{-p})}{p}=2(1-e^{-p})\le 2p.  
\end{equation}
Using this coupling, the vertices of $\bM_p$ form a subset of the
vertices of $\bX_p$, possibly with some repetitions,  and it follows
that if vertex $i$ is selected for $\bX_p$ but is isolated, and thus does
not appear in $X_p$, then it also will not appear in $M_p$, since even if it
appears in $\bM_p$ it will be isolated there. Hence, if $X_p\neq M_p$, there
must be some vertex $i\in X_p$ such that $Y_i\neq I_i$.
Consequently, by \eqref{nora},
\begin{equation}
  \P\bigpar{X_p\neq M_p\mid X_p}
\le \sum_{i\in V(X_p)} \P\bigpar{Y_i\neq I_i\mid I_i=1}
\le 2p v(X_p).
\end{equation}
\end{proof}

As in \cite{VR2}, we obtain a corollary on stretched convergence of
unstretched empirical graphons, \cf{} \refSS{SSstretch}.
\begin{theorem}[{Extension of \cite[Theorem 5.7(2)]{VR2}}]\label{T3}
If\/ $\cW$ is any non-zero graphex, then
$\hW_{G_{\tau_k}(\cW),1}   \togs \cW$.
\end{theorem}
\begin{proof}
  Immediate by \refT{T2} and \refL{Ltogs}.
\end{proof}

\begin{remark}
Note that even if $\cW=(I,S,W)$ is a general graphex with non-zero dust and star
components $I$ and $S$, 
\refT{T2} and \cite[Theorem 4.12]{VR2} yield convergence to $\cW$ of a
sequence of graphons, where the dust and star components are taken to be 0.
This shows that 
for finite $r$, it is not possible to distinguish between isolated edges or
stars coming from $I$ and $S$ and isolated edges or stars coming from the
graphon part $W$.
Note that in contrast,
for the infinite random graph $\gG(\cW)$, there is
an obvious difference between edges produced by $I$, $S$ and $W$: they have
\as{} 2, 1 and 0 endpoints of degree 1, respectively.
\end{remark}

\section*{Acknowledgement}
I thank Daniel Roy for an interesting discussion
at the Isaac Newton Institute for Mathematical Sciences
during the programme 
Theoretical Foundations for Statistical Network Analysis in 2016
(EPSCR Grant Number EP/K032208/1).
This work was also supported by the Knut and Alice Wallenberg Foundation.

\newcommand\AAP{\emph{Adv. Appl. Probab.} }
\newcommand\JAP{\emph{J. Appl. Probab.} }
\newcommand\JAMS{\emph{J. \AMS} }
\newcommand\MAMS{\emph{Memoirs \AMS} }
\newcommand\PAMS{\emph{Proc. \AMS} }
\newcommand\TAMS{\emph{Trans. \AMS} }
\newcommand\AnnMS{\emph{Ann. Math. Statist.} }
\newcommand\AnnPr{\emph{Ann. Probab.} }
\newcommand\CPC{\emph{Combin. Probab. Comput.} }
\newcommand\JMAA{\emph{J. Math. Anal. Appl.} }
\newcommand\RSA{\emph{Random Struct. Alg.} }
\newcommand\ZW{\emph{Z. Wahrsch. Verw. Gebiete} }
\newcommand\DMTCS{\jour{Discr. Math. Theor. Comput. Sci.} }

\newcommand\AMS{Amer. Math. Soc.}
\newcommand\Springer{Springer-Verlag}
\newcommand\Wiley{Wiley}

\newcommand\vol{\textbf}
\newcommand\jour{\emph}
\newcommand\book{\emph}
\newcommand\inbook{\emph}
\def\no#1#2,{\unskip#2, no. #1,} %(typeset after year) 
\newcommand\toappear{\unskip, to appear}

\newcommand\arxiv[1]{\texttt{arXiv:#1}}
\newcommand\arXiv{\arxiv}

\def\nobibitem#1\par{}

\end{document}